\documentclass[11pt]{article}

\usepackage{lipsum}

\newcommand\blfootnote[1]{%
  \begingroup
  \renewcommand\thefootnote{}\footnote{#1}%
  \endgroup
}

\usepackage[utf8]{inputenc} %

\usepackage{fancyhdr} 
\usepackage{amsthm}

\usepackage[utf8]{inputenc} 
\usepackage[T1]{fontenc}
\usepackage{url}
\usepackage{ifthen}

\usepackage{algpseudocode,algorithm,algorithmicx}

\usepackage{graphicx}
\usepackage{amssymb}
\usepackage{xcolor}
\usepackage{enumitem}
\usepackage{mathtools}
\usepackage{latexsym}
\usepackage{amsmath}
\usepackage{amsfonts}
\usepackage{amsthm}
\usepackage{amssymb}
\usepackage{bbm}
\usepackage{centernot}

\usepackage{listings}
\usepackage{color} %
\definecolor{mygreen}{RGB}{28,172,0} %
\definecolor{mylilas}{RGB}{170,55,241}
\usepackage{thm-restate}

\usepackage[colorlinks,citecolor=blue,linkcolor=magenta,bookmarks=true,hypertexnames=false]{hyperref}
\usepackage[nameinlink,capitalise]{cleveref}

\crefname{equation}{Equation}{Equations}
\crefname{lemma}{Lemma}{Lemmata}
\crefname{claim}{Claim}{Claims}
\crefname{theorem}{Theorem}{Theorems}
\crefname{proposition}{Proposition}{Propositions}
\crefname{corollary}{Corollary}{Corollaries}
\crefname{claim}{Claim}{Claims}
\crefname{remark}{Remark}{Remarks}
\crefname{definition}{Definition}{Definitions}
\crefname{fact}{Fact}{Facts}
\crefname{question}{Question}{Questions}
\crefname{condition}{Condition}{Conditions}
\crefname{algorithm}{Algorithm}{Algorithms}
\crefname{assumption}{Assumption}{Assumptions}

\newtheorem{theorem}{Theorem}[section]

\newtheorem{fact}[theorem]{Fact}

\theoremstyle{definition}

\newtheorem{remark}[theorem]{Remark}

\newcommand{\eps}{\epsilon}

\newcommand{\Var}{\mathbf{Var}}

\def\P{\mathbb P}
\def\R{\mathbb R}

\def\sgn{\mathrm{sgn}}

\DeclareMathOperator*{\E}{\mathbf{E}}

\newcommand{\cN}{\mathcal{N}}

\def\colorful{0}

\ifnum\colorful=1
\newcommand{\new}[1]{{\color{red} #1}}

\else
\newcommand{\new}[1]{{#1}}

\fi

\newcommand{\GMT}{Gaussian-Mean-Testing}

\allowdisplaybreaks

\title{Gaussian Mean Testing Made Simple\blfootnote{Authors are in alphabetical order.}}

\author{
Ilias Diakonikolas\thanks{Supported by NSF Medium Award CCF-2107079,
NSF Award CCF-1652862 (CAREER), a Sloan Research Fellowship, and
a DARPA Learning with Less Labels (LwLL) grant.}\\
University of Wisconsin-Madison\\
{\tt ilias@cs.wisc.edu}\\
\and
Daniel M. Kane\thanks{Supported by NSF Medium Award CCF-2107547,
NSF Award CCF-1553288 (CAREER), and a Sloan Research Fellowship.}\\
University of California, San Diego\\
{\tt dakane@cs.ucsd.edu}
\and
Ankit Pensia\thanks{Supported by NSF Award CCF-1652862 (CAREER), 
and NSF grants CCF-1841190 and CCF-2011255.}\\
University of Wisconsin-Madison\\
{\tt ankitp@cs.wisc.edu}
}

\begin{document}

\maketitle

\setcounter{page}{0}

\thispagestyle{empty}

\begin{abstract}
We study the following fundamental hypothesis testing problem, which we term {\em Gaussian mean testing}. 
Given i.i.d.\ samples from a distribution $p$ on $\R^d$, the task is to distinguish, with high probability, 
between the following cases:
(i) $p$ is the standard Gaussian distribution, $\cN(0,I_d)$, and
(ii) $p$ is a Gaussian $\cN(\mu,\Sigma)$ for some unknown covariance $\Sigma$ 
and mean $\mu \in \R^d$ satisfying $\|\mu\|_2 \geq \epsilon$.
Recent work gave an algorithm for this testing problem with the 
optimal sample complexity of $\Theta(\sqrt{d}/\epsilon^2)$.
Both the previous algorithm and its analysis are quite complicated.
Here we give an extremely simple algorithm for Gaussian mean testing 
with a one-page analysis. Our algorithm is sample optimal
and runs in sample linear time.
\end{abstract}

\newpage

\section{Introduction} %
\label{sec:introduction}

The paradigmatic problem in distribution testing~\cite{GRexp:00, BFR+:00} 
is the following: given sample access to an unknown distribution
$p$, determine whether $p$ has some global property or is ``far''
from any distribution having the property.
During the past two decades, a wide range of properties have been studied,
and we now have sample-optimal testers 
for many of them~\cite{Paninski:08, CDVV14, VV14, ADK15, DK:16, 
CDGR18, DiakonikolasGPP18, CanonneDKS18, DiakonikolasGKP21}.

Without a priori assumptions on the underlying distribution $p$,
at least $\Omega(\sqrt{N})$ many samples are required 
for testing even the simplest properties, where $N$ is the domain size of $p$.
If $p$ is either high-dimensional 
(supported on an exponentially large domain, e.g., $\{0, 1\}^d$) 
or continuous, such a sample bound is prohibitive.
This observation has motivated a line of work studying 
distribution testing of {\em structured} distribution families. 
This includes both nonparametric families in 
low-dimensions~\cite{DKN:15, DKN:15:FOCS, DKN17, DiakonikolasKP19}
and parametric families in high dimensions~\cite{CDKS17, DasDK19, AcharyaBDK18, CCKLW21-subcube}.

This work focuses on the high-dimensional setting.
Arguably the most basic high-dimensional testing problem is the following: We assume that
$p$ is an {\em identity covariance} Gaussian distribution on $\R^d$ and the goal is to distinguish between the
cases that its mean is zero or at least $\eps$ in $\ell_2$-norm. This is known as the {\em Gaussian sequence model} 
and has a rich history in statistics~\cite{Ermakov91,Baraud02,IngSus03} 
(see also~\cite{DKS17-sq} for a simple algorithm and matching lower bound).
In particular, the aforementioned works have established that the sample complexity 
of this basic problem is $\Theta(\sqrt{d}/\epsilon^2)$.

Here we consider a generalization of the Gaussian sequence model, 
recently studied in the TCS literature~\cite{CCKLW21-subcube}. 
We will call this problem \GMT:
\vspace{0.2cm}

\fbox{\parbox{6in}{
{\bf Problem}: \GMT\\

\vspace{-0.2cm}

\textbf{Input:} Sample access to a distribution $p$ supported on $\R^d$ and $\epsilon>0.$
\vspace{0.1cm}

\textbf{Output:} 
\begin{itemize}
   \item \new{(Completeness)} ``ACCEPT'' with probability at least $2/3$ if $p = \cN(0,I_d)$,
 
  \item \new{(Soundness)} ``REJECT'' with probability at least $2/3$ if $p = \cN(\mu,\Sigma)$ for some $\mu$ with $\|\mu\|_2 \geq \epsilon$. 
\end{itemize}
}}
\vspace{0.2cm}

The work of \cite{CCKLW21-subcube} proposed a testing algorithm for this problem 
with the optimal sample complexity of $n = \Theta(\sqrt{d}/\epsilon^2)$ 
for the (restricted) parameter regime of $\eps \in (0,1]$. The initial algorithm of \cite{CCKLW21-subcube}
had quasi-polynomial time complexity; this was improved to polynomial using 
an observation from a subsequent work \cite{CJLW21-junta}.

Before we describe our results, we summarize the approach of \cite{CCKLW21-subcube}.

\paragraph{Approach of \cite{CCKLW21-subcube}} 
The algorithm of \cite{CCKLW21-subcube} proceeds by reducing {\GMT} 
to another high-dimensional testing problem, which we call Hypercube-Mean-Testing, described next.
Given sample access to a distribution $p$ on the hypercube $\{ \pm 1\}^d$ 
and $\eps \in (0,1)$, Hypercube-Mean-Testing asks us to distinguish 
between the cases when $p$ is uniform and 
when the mean of the distribution $p$ has $\ell_2$-norm at least $\epsilon$.  

The relation between Hypercube-Mean-Testing and \GMT\ is apparent 
when we look at the function $F: \R^d \to \{\pm 1\}^d$ 
that maps $x \in \R^d$ to $y \in \{\pm 1\}^d$ coordinatewise by $y_i = \sgn(x_i)$. 
For a distribution $p$, we use $F(p)$ to denote the distribution of $F(X)$ when $X \sim p$.
It is not too hard to see that when $p = \cN(0,I_d)$, then $F(p)$ is the uniform distribution on the hypercube. 
Similarly, it can be shown that when $p \sim \cN(\mu, \Sigma)$ 
with $\|\mu\|_2 \geq \eps$ for $\epsilon \in (0,1]$ and $\max_{i \in [d]} \Sigma_{i,i} \leq 2$, then 
$F(p)$ is a distribution on hypercube whose mean has euclidean norm larger than $\Omega(\eps)$.  
Thus the algorithm of \cite{CCKLW21-subcube} works as follows:
  \begin{enumerate}
     \item Check if there is a coordinate $i \in [d]$ such that $\Sigma_{i,i} \geq 2$, which can be tested reliably and efficiently with $O(\log d)$ samples, and
     \item Run a tester for Hypercube-Mean-Testing with input distribution $F(p)$ and $\eps' = \Theta(\eps)$.
   \end{enumerate} 
One of the main contributions of \cite{CCKLW21-subcube} is an algorithm for Hypercube-Mean-Testing 
with sample complexity $n = O(\sqrt{d}/\eps^2)$ 
and runtime\footnote{This follows by the reduction procedure of  \cite{CCKLW21-subcube} with a more careful runtime analysis by \cite{CJLW21-junta}.} $O(n^2d)$. 
The proposed algorithm and its analysis are somewhat involved. 
Moreover, the reduction of \GMT\ to Hypercube-Mean-Testing 
crucially uses that the parameter $\epsilon$ in {\GMT} 
is less than a small enough constant.
Our main contribution is a very simple direct algorithm
with a compact analysis 
(by avoiding this reduction, our algorithm works for all values of $\epsilon$).

\paragraph{Our Result}
In this paper, we give a very simple sample-optimal algorithm (\Cref{alg:tester}) for \GMT \
whose full analysis fits in one page. In addition to its simplicity, our algorithm
runs in linear time (requires a single pass over the data) and works for all values of $\epsilon$. 

In particular, we establish the following result:
\begin{theorem}
\label{thm:GaussMeanTesting}
\Cref{alg:tester} solves {\GMT} with $n = \Theta\left( \max(1, \sqrt{d}/\epsilon^2)\right)$ samples 
and can be implemented in $O(n \, d)$ time.
\end{theorem}

Our algorithm is extremely easy to describe: 
we sample two sets of $\Theta\left( \max(1, \sqrt{d}/\epsilon^2)\right)$ samples, 
compute the sample mean of each of them, 
and calculate the inner product of the two obtained vectors. 
The algorithm outputs ``ACCEPT'' if the inner product 
has small absolute value, and ``REJECT'' otherwise. 
A detailed pseudocode follows.

 \begin{algorithm}[h]  
    \caption{GaussianMeanTester} 
      \label{alg:tester}
    \begin{algorithmic}[1]
      \Statex \textbf{Input}: 
      Sample access to distribution $p$ on $\R^d$ and $\epsilon > 0$.
      \Statex \textbf{Output}: 
      ``ACCEPT'' if $p = \cN(0,I_d)$, ``REJECT'' if $p = \cN(\mu,\Sigma)$ and $\|\mu\|_2 \geq \epsilon$; both with probability at least $2/3$.

      \State Set $n = 25 C_*^2 \sqrt{d} / \epsilon^2 $, where $C_*$ is the absolute constant from \Cref{ThmCarWri}.
      \State Sample $2n$ i.i.d.\ points from $p$ and denote them by $X_1,\dots,X_n$ and $Y_1,\dots,Y_n$.
      \State Define $Z = (1/n^2) (\sum_{i=1}^n X_i)^\top (\sum_{i=1}^n Y_i)$.
       \If{$|Z| \leq \sqrt{3d}/n$}
       \State \Return ``ACCEPT''
       \Else 
       \State \Return ``REJECT''
       \EndIf   

    \end{algorithmic}  
  \end{algorithm}

We now describe the high-level idea of our proof.
\paragraph{Our Technique}
A natural first attempt at a tester would be to attempt 
to approximate $\|\mu\|_2^2$ and to reject if the answer is too large.
A reasonable way to do this is to compute two independent estimates 
of $\mu$ and take their inner product. 
For example, as $Z = (\sum_{i=1}^n X_i/n)^\top (\sum_{i=1}^n Y_i/n)$, 
where $X_1, \dots, X_n$ and $Y_1, \dots, Y_n$ are all independent samples.
(Note that the two estimates above should be independent: 
one cannot take $(\sum_{i=1}^n X_i/n)^\top (\sum_{i=1}^n X_i/n)$ 
as a reliable estimate of $\|\mu\|_2^2$ because 
the correlation between the two halves will bias the final estimate.)
Computing the variance of $Z$, $\Var[Z]$, it is not hard to see that 
in the completeness case, $|Z| =O(\sqrt{d}/n)$ with high probability.
We would like to be able to claim that in the soundness case the quantity $|Z - \|\mu\|_2^2|$ is similarly small;
this would give us an easy separation, 
so long as $\sqrt{d}/n \leq c \cdot \epsilon^2$ 
for a small enough constant $c > 0$.
Unfortunately, this may not hold if $\Sigma$ is too large, 
as this may cause the variance of $Z$ to be larger than required.
Fortunately, in this case we are rescued by another argument. 
If $\Var[Z]$ is large, it means that the values of $Z$ will be spread over a very large range.
This in turn will make it unlikely that $Z$ will be in the narrow range of values where $|Z| = O(\sqrt{d}/n)$.
In particular, we can formalize this by noting that $Z$ is a degree-$2$ polynomial in Gaussian inputs, 
and applying the Carberry-Wright anti-concentration inequality (see \cref{ThmCarWri}). 
Thus, we end up with the simple tester of compute $Z$, 
and then check whether or not $|Z| < C \sqrt{d}/n$ 
for an appropriate constant $C$.

\new{ \begin{remark} Our testing algorithm applies 
for the following more general testing problem: 
\begin{enumerate}[leftmargin=*]
  \item (Completeness) $p$ is a distribution with  mean $\mu$ and covariance $\Sigma$ 
  satisfying $\|\mu\|_2 \leq c\cdot\epsilon$ for a small enough positive constant $c<1$ and $\|\Sigma\|_{\mathrm{F}} \leq \sqrt{d}$, 
  where $\|\cdot\|_{\mathrm{F}}$ is the Frobenius norm.
  \item (Soundness) $p$ is a log-concave distribution with mean $\mu$ satisfying $\|\mu\|_2 \geq \epsilon$.
\end{enumerate}   
\end{remark}
}
\section{Proof of \cref{thm:GaussMeanTesting}}

\paragraph{Notation and Background}
We use $I_d$ to denote the $d\times d$ identity matrix.
For a vector $x \in \R^d$, we use $\|x\|_2$ to denote its Euclidean norm.
For a univariate random variable $X$, we use $\E[X]$ and $\Var[X]$ 
to denote its mean and variance, respectively. 
The multivariate Gaussian distribution with mean $\mu$ and covariance $\Sigma$ 
is denoted by $\cN(\mu, \Sigma)$.
For two matrices $A$ and $B$, we use $\langle A, B\rangle$ 
to denote the trace inner product, i.e., $\langle A, B\rangle = \mathrm{tr}(A^\top B)$.

We will require the following well-known fact.

\begin{fact}[Carbery-Wright inequality for quadratics~\cite{CarWri01}]
\label{ThmCarWri}
There exists a $C_{*} > 0$ such that the following holds:  
Let $G \sim \cN(0,I_d)$ in $\R^d$,  $p: \R^{d} \to \R$ be a degree-$2$ polynomial, 
and $\alpha \in (0,\infty)$. Then we have that
$ \P(|p(G)| \leq \alpha \sqrt{\E [p^2(G)]}) \leq C_{*}  \sqrt{\alpha}$. 
\end{fact}

We are now ready to prove our main result.
\begin{proof}(of \cref{thm:GaussMeanTesting})
\cref{alg:tester} samples $2n$ points, where $n \geq 25C_*^2 \sqrt{d} / \epsilon^2$ and $C_*$ 
is the constant in the Carbery-Wright Theorem (cf. \cref{ThmCarWri});
 without loss of generality, we will assume that $n \geq 1$.  
Let the $2n$ samples be $X_1,\dots,X_n$ and $Y_1,\dots,Y_n$, 
where each $X_i$ and $Y_i$ is distributed as  $\cN(\mu,\Sigma)$.
We define our test statistic to be $Z := (1/n^2) (\sum_{i=1}^n X_i)^\top(\sum_{i=1}^n Y_i)$.
\cref{alg:tester}  outputs ``ACCEPT'' if $|Z| \leq \sqrt{3d}/n$ and outputs 
``REJECT'' otherwise. 
The claim on the running time is immediate.  
We now analyze the  completeness (the algorithm outputs ``ACCEPT'' with probability $2/3$ when $p = \cN(0,I_d)$)
and soundness  (the algorithm outputs ``REJECT'' with probability $2/3$ 
when $p = \cN(\mu,\Sigma)$ and $\|\mu\|_2 > \epsilon$) of the proposed test.

\paragraph*{Completeness}

We begin by calculating the mean and variance of $Z$ when $p = \cN(0,I_d)$. 
Let $G_1$ and $G_2$ be two independent $\cN(0,I_d)$ random variables.
Since $(\sum_i X_i)/\sqrt{n}$ and $(\sum_i Y_i)/\sqrt{n}$ 
are independently distributed as $\cN(0, I_d)$, 
it follows that $Z$ has the same distribution as $(1/n) (G_1 ^\top G_2)$.
This directly gives use that $\E [Z] = 0$.
The variance of $Z$ can be calculated as follows:
\begin{align*}
 \Var[Z]  &= 	\frac{1}{n^2}\E \left[\left(  G_1^\top G_2 \right)^2\right] = \frac{1}{n^2} \E \left[\left\langle  G_1G_1^\top,  G_2G_2^\top \right\rangle \right] 	= \frac{1}{n^2} \langle I_d, I_d \rangle	=    \frac{d}{n^2},
 \end{align*}
 where we use that $G_1$ and $G_2$ are independent.
By Chebyshev's inequality, with probability at least $2/3$, $|Z| \leq \sqrt{3 d}/n$. 
Therefore, the probability of acceptance when $p = \cN(0,I_d)$ is at least $2/3$.

\paragraph{Soundness}
We now consider the setting when $\|\mu\|_2 \geq \epsilon$ and $\Sigma$ is an arbitrary positive semidefinite matrix. 
Since $X_i$'s and $Y_i$'s are i.i.d.\ with mean $\mu$, 
we get that $\E[Z] = \|\mu\|_2^2$, which is larger than $\epsilon^2$.
Our goal will be to show that, with probability at least $2/3$, $|Z| > 2 \sqrt{d}/n$.
We will use the Carbery-Wright inequality (\Cref{ThmCarWri}).
Observe that $Z$ is a quadratic polynomial of a Gaussian distribution.
Let $\|Z\|_{L_2}$ denote $\sqrt{\E [Z^2]}$ and observe that $\|Z\|_{L_2} \geq \E[Z] \geq \epsilon^2$. 
By \Cref{ThmCarWri}, we obtain the following:
\begin{align*}
\P\left(|Z| \leq 2\sqrt{d}/n\right) &= \P\left(|Z| \leq  \frac{2\sqrt{d}/n}{\|Z\|_{L_2}} \cdot \|Z\|_{L_2}\right) \leq C_* \sqrt{\frac{2 \sqrt{d}/n}{\|Z\|_{L_2}}}  \leq 	C_*\sqrt{\frac{ 2 \sqrt{d}/n}{ \epsilon^2}}  < \frac{1}{3}, 
\end{align*}
where the second inequality uses that $\|Z\|_{L_2} \geq \epsilon^2$ and the last inequality uses that $n \geq 25C_*^2 \sqrt{d}/\epsilon^2$.
Therefore, the probability of rejection, i.e., of the event where $|Z| > \sqrt{3 d}/n$, is at least $2/3$.
Combining the above, it follows that the algorithm correctly rejects with probability at least $2/3$.
\end{proof}

\bibliographystyle{alpha}
\bibliography{allrefs}

\appendix

\end{document}